\documentclass[11pt,reqno]{amsart}
\usepackage{graphicx}
\usepackage{verbatim}
\usepackage{textcomp}
\usepackage{amssymb}
\usepackage{cite}
\usepackage{hyperref}
\usepackage{amsmath}
\usepackage{latexsym}
\usepackage{amscd}
\usepackage{amsthm}
\usepackage{mathrsfs}
\usepackage{xypic}
\usepackage{bm}
\usepackage{url}

\newtheorem{thm}{Theorem}[section]
\newtheorem{corr}[thm]{Corollary}
\newtheorem{lem}[thm]{Lemma}
\newtheorem{prop}[thm]{Proposition}

\newtheorem{exam}{Example}[section]
\theoremstyle{definition}

\theoremstyle{remark}
\newtheorem{rem}{Remark}[section]
\numberwithin{equation}{section}
\def\R{\mathbb R}

\def\S{\mathbb S}
\def\cal{\mathcal}
\def\f{\frac}
\def\td{\tilde}
\def\ra{\rightarrow}

\begin{document}
\title[Classification and rigidity of Self-shrinkers]{Classification and rigidity of self-shrinkers in the Mean curvature flow}
\author{Haizhong Li}
\address{Department of mathematical sciences, and Mathematical Sciences
Center, Tsinghua University, 100084, Beijing, P. R. China}
\email{hli@math.tsinghua.edu.cn}
\author{Yong Wei}
\address{Department of Mathematical Sciences, Tsinghua University, 100084, Beijing, China.}
\email{wei-y09@mails.tsinghua.edu.cn}

\subjclass[2010]{Primary 53C42; Secondary 53C44}


\keywords {Self-shrinker, rigidity, mean curvature flow}

\thanks{The authors were supported by NSFC No.11271214 and Tsinghua University-K. U. Leuven Bilateral Scientific Cooperation Fund.}
\maketitle

\begin{abstract}
In this paper, we first use the method of Colding and Minicozzi II \cite{CM2} to show that K. Smoczyk's classification theorem \cite{Smo}
for complete  self-shrinkers in higher codimension also holds under a
weaker condition. Then as an application, we give some rigidity results for self-shrinkers in arbitrary codimension.
\end{abstract}

\section {Introduction}

An immersion $x:M^n\to \R^{n+p}$ of a smooth n-dimensioal manifold $M$ into the Euclidean space is called a
self-shrinker if it satisfies the quasilinear elliptic system:
\begin{equation}\label{1-1}
{\bf H}=-x^{\perp},
\end{equation}
where ${\bf H}$ denotes the mean curvature vector of the immersion
and $\perp$ is the projection onto the normal bundle of $M$.

Self-shrinkers play an important role in the study of  the mean
curvature flow. Not only they correspond to self-shrinking solutions to the mean curvature flow, but also they describe all possible blow ups at a given singularity of the mean curvature flow. We refer the readers to \cite{CM2,H2,H3,I,White} and references therein for more information on self-shrinkers and singularities of mean curvature flow.

There are many results about the classification of self-shrinkers. In the curve case, U. Abresch and J. Langer \cite{AbL} gave a complete
classification of all solutions to \eqref{1-1}. These curves are now
called Abresch-Langer curves, and the only simple closed one is the circle.

In higher dimension and codimension one, Huisken \cite{H2,H3} (see also \cite{Zhu}) proved a classification theorem
for smooth mean convex self-shrinkers $M^n$ in $\R^{n+1}$ with polynomial
volume growth. Suppose further in the noncompact case $|(\nabla)^kA|$ are uniformly bounded for $k=0,1,2$, as well as $|A|^2\leq CH^2$ everywhere on $M$.
Then $M$ are isometric to $\Gamma\times \R^{n-1}$ or
$\S^k(\sqrt{k})\times\R^{n-k}$ ($0<
 k\leq n$). Here, $\Gamma$ is a Abresch-Langer curve and $\S^k(\sqrt{k})$ is a
$k$-dimensional sphere. Recently, Colding and Minicozzi II
\cite{CM2} generalized this result and showed that Huisken's classification theorem still
holds without the assumption on the bounds for derivatives of the second fundamental form. Moreover,
they showed that the only smooth embedded entropy stable
self-shrinkers with polynomial volume growth in $\R^{n+1}$ are
hyperplanes, n-spheres, and cylinders. By imposing symmetries, Kleene and M{\o}ller \cite{KM} classified the complete $n$-dimensional embedded self-shrinkers of revolution in $\R^{n+1}$.

In arbitrary codimension, the situation becomes more complicated. K. Smoczyk  \cite{Smo} proved the following results:

\medskip\noindent
{\bf Theorem A}(Theorem 1.1 in \cite{Smo}).  {\it
Let $x:M^n\ra \R^{n+p}$ be a closed self-shrinker, then $M$ is a minimal submanifold of the sphere $\S^{n+p-1}(\sqrt{n})$ if and only if ${\bf H}\neq 0$ and $\nabla^{\perp}\nu=0$, where $\nu={\bf H}/|{\bf H}|$ is the principal normal.
}

\medskip\noindent
{\bf Theorem B}(Theorem 1.2 in \cite{Smo}).  {\it
Let $x:M^n\ra \R^{n+p}$ be a complete non-compact connected self-shrinker with ${\bf H}\neq 0$ and $\nabla^{\perp}\nu=0$. Suppose further that $M$ has uniformly bounded geometry, that is, there exists constants $c_k$ such that $|(\nabla)^k A|\leq c_k$ for any $k\geq 0$. Then $M$ must belong to one of the followings:
\begin{equation*}
    \Gamma\times \R^{n-1},\qquad \td{M}^r\times\R^{n-r}.
\end{equation*}
Here, $\Gamma$ is one of the Abresch-Langer curves and $\td{M}^r$ is a complete minimal submanifold of the sphere $\S^{r+p-1}(\sqrt{r})\subset\R^{p+r}$, where $0<r=\textrm{rank}(A^{\nu})\leq n$ denotes the rank of the principal second fundamental form $A^{\nu}=<\nu,A>$.
}

\medskip
 Note that any blowup of a Type-I singularity of the mean curvature flow forming on a compact submanifold will automatically be complete with uniformly bounded geometry, therefore Theorem B may be applied to those blowup limits. In Smoczyk's proof of Theorem B, the uniformly bounded geometry is needed in the integrating by parts with respect to the Gauss kernel $\rho(x)=e^{-\f 12|x|^2}$. In the first part of this paper, we will use the method of Colding and Minicozzi II \cite{CM2} to show that Theorem B also holds under a weaker condition:

\begin{thm}
Let $x:M^n\ra \R^{n+p}$ be a complete non-compact connected self-shrinker with ${\bf H}\neq 0$ and $\nabla^{\perp}\nu=0$. Suppose further that $M$ has polynomial volume growth and satisfies $|A|^2-|A^{\nu}|^2\leq c$ for some constant c, where $A^{\nu}=<\nu,A>$ is the principal second fundamental form. Then $M$ must belong to one of the followings:
\begin{equation}
    \Gamma\times \R^{n-1},\qquad \td{M}^r\times\R^{n-r}.
\end{equation}
Here, $\Gamma$ is one of the Abresch-Langer curves and $\td{M}^r$ is a complete minimal submanifold of the sphere $\S^{p+r-1}(\sqrt{r})\subset\R^{p+r}$, where $0<r=\textrm{rank}(A^{\nu})\leq n$ denotes the rank of $A^{\nu}$.
\end{thm}

\begin{rem}
In the recent paper \cite{ALW}, Ben Andrews and the authors considered the $\mathcal{F}$-stability of self-shrinkers in arbitrary codimension, where Theorem A and Theorem 1.1 applied.
\end{rem}

When $p=1$ and the self-shrinker is embedded, Theorem A and Theorem 1.1 reduce to Colding-Minicozzi II's result (Theorem 0.17 in \cite{CM2}):

\begin{corr}
$\S^k(\sqrt{k})\times \R^{n-k}$ are the only complete embedded self-shrinkers without boundary, with polynomial volume growth, and $H\geq 0$ in $\R^{n+1}$.
\end{corr}

In the statement of Theorem 1.1, We say that a submanifold $M^n$ in $\R^{n+p}$  has polynomial volume growth if there exists constants $C$ and $d$ such that for all $r\geq 1$, there holds
\begin{equation*}
    \textrm{Vol} (B(r)\cap M)\leq Cr^d,
\end{equation*}
where $B_r$ denotes an Euclidean ball with radius $r$. By using Huisken's monotonicity formula \cite{H2}, Colding-Minicozzi II \cite{CM2} proved that any self-shrinker which arises as the blow up at a given singularity in the mean curvature flow must have polynomial volume growth. Their result is proved for self-shrinkers in the hypersurface case, but it also holds for arbitrary codimension.

In the second part of this paper, we shall apply Theorem 1.1 to give some rigidity properties of self-shrinkers with higher codimension. Recall that the first gap theorem for self-shrinkers was proved by N. Q. Le and N. Sesum \cite{LeS} for hypersurface case. Later this was generalized by H. -D. Cao and H. Li \cite{CL2} to self-shrinkers with arbitrary codimension, they showed that

\medskip\noindent
{\bf Theorem C.}( Cao-Li \cite{CL2}) {\it Let $x: M^n\ra\R^{n+p}$ be an n-dimensional complete self-shrinker without
boundary and with polynomial volume growth, if
\begin{equation}
   0\leq |A|^2\leq 1,
\end{equation}
then either (i) $|A|^2\equiv 0$ and $M$ is the hyperplane, or (ii) $|A|^2\equiv 1$ and $M$ is $\S^m(\sqrt{m})\times\R^{n-m}$ in $\R^{n+1}$ with $1\leq m\leq n$.
 }

\vskip 2mm
\begin{rem}
We remark that Theorem C is independent of the dimension and codimension of the self-shrinker. In \cite{Cheng-P}, the authors got a related result without condition ``polynomial volume growth".
\end{rem}

Although classifying the self-shrinker with higher codimension is complicated,  we can also apply Theorem A and Theorem 1.1 to give some results for self-shrinkers in some special situations. In the following we will consider the self-shrinkers with codimension 2, 2-dimension and with flat normal bundle respectively.

\begin{thm}
Let $X:M^n\ra\R^{n+2}$ be a complete embedded self-shrinker without boundary and with polynomial volume growth. If ${\bf H}\neq 0$, $\nabla^{\perp}\nu=0$, where $\nu$ is the principal normal, and
\begin{equation}
    1\leq|A|^2\leq 2,
\end{equation}
then there are two possibilities:

\vskip 1mm

(i) $|A|^2\equiv 1$ and $M$ is $\S^m(\sqrt{m})\times\R^{n-m}$ in $\R^{n+1}$ with $1\leq m\leq n$.

\vskip 1mm

(ii) $|A|^2\equiv 2$ and $M$ is one of the self-shrinkers in Example 1.1 below.
\end{thm}

\begin{thm}
Let $X:M^n\ra\R^{n+2}$ be a complete embedded self-shrinker without boundary and with polynomial volume growth. If ${\bf H}\neq 0$, $\nabla^{\perp}\nu=0$, where $\nu$ is the principal normal. Then there exists a constant $\delta>0$ such that if
\begin{equation}\label{th1-3}
    2\leq|A|^2\leq 2+\delta,
\end{equation}
then $|A|^2\equiv 2$ and $M$ is one of  the self-shrinkers in Example 1.1 below.
\end{thm}

\begin{exam}
Let
\begin{equation*}
    \td{M}^r=\S^k(\sqrt{k})\times\S^{r-k}(\sqrt{r-k})\hookrightarrow\S^{r+1}(\sqrt{r}),\quad 1\leq k\leq r-1
\end{equation*}
be the Clifford minimal hypersurfaces in the sphere $\S^{r+1}(\sqrt{r})$. Then
\begin{equation*}
    x:M^n=\td{M}^r\times\R^{n-r}\hookrightarrow\R^{n+2},\qquad (2\leq r\leq n)
\end{equation*}
is a complete embedded self-shrinker without boundary and with polynomial volume growth, with parallel principal normal and $|A|^2\equiv 2$.
\end{exam}

\begin{exam}
Let
\begin{equation*}
    \td{M}^r\hookrightarrow\S^{r+1}(\sqrt{r})
\end{equation*}
be an isoparametric minimal hypersurface, then
\begin{equation*}
    x:M^n=\td{M}^r\times\R^{n-r}\hookrightarrow\R^{n+2},\qquad (1\leq r\leq n)
\end{equation*}
is a self-shrinker with $|A|^2$ can only be $1,2,3,4,6$.
\end{exam}

\vskip 2mm
For the 2-dimensional self-shrinkers, we have the following two rigidity results.
\begin{thm}
Let $X:M^2\ra\R^{2+p}$ be a 2-dimensional complete embedded self-shrinker without boundary and with polynomial volume growth. If ${\bf H}\neq 0$, $\nabla^{\perp}\nu=0$, where $\nu$ is the principal normal, and
\begin{equation}\label{dim2}
   1\leq|A|^2\leq \f 53,
\end{equation}
then there are two possibilities:

\vskip 1mm

(i) $|A|^2\equiv 1$ and $M=\S^2(\sqrt{2})$  or $\S^1(1)\times\R$ in $\R^3$.

\vskip 1mm

(ii) $|A|^2\equiv \f 53$ and $M$ is the self-shrinker in Example 1.3 below.
\end{thm}

\begin{thm}
Let $X:M^2\ra\R^{2+p}$ be a 2-dimensional complete embedded self-shrinker without boundary and with polynomial volume growth. If ${\bf H}\neq 0$, $\nabla^{\perp}\nu=0$, where $\nu$ is the principal normal, and
\begin{equation}\label{dim2-1}
   \f 53\leq|A|^2\leq \f {11}6,
\end{equation}
then there are two possibilities:

\vskip 1mm

(i) $|A|^2\equiv  \f 53$ and $M$ is the self-shrinker in Example 1.3 below.

\vskip 1mm

(ii) $|A|^2\equiv \f {11}6$ and $M$ is the self-shrinker in Example 1.4 below.
\end{thm}

\begin{exam}
The canonical minimal immersion (see \cite{KS}, \cite{LS})
\begin{equation*}
    x:\S^{2}(\sqrt{6})\ra\S^4(\sqrt{2})
\end{equation*}
has $|\td{A}|^2\equiv \f 23$. $x(\S^{2}(\sqrt{6}))\subset\S^4(\sqrt{2})$ is called the Veronse surface. Consider it as a submanifold in $\R^5$,
\begin{equation*}
    x:\S^{2}(\sqrt{6})\ra\S^4(\sqrt{2})\hookrightarrow\R^5
\end{equation*}
it is a self-shinker with $|A|^2\equiv 1+\f 23=\f 53$.
\end{exam}

\begin{exam}
The canonical minimal immersion (see \cite{KS})
\begin{equation*}
    x:\S^2(\sqrt{12})\ra\S^6(\sqrt{2})
\end{equation*}
has $|\td{A}|^2\equiv \f 56$, consider it as a submanifold in $\R^7$
\begin{equation*}
    x:\S^2(\sqrt{12})\ra\S^6(\sqrt{2})\hookrightarrow\R^7
\end{equation*}
it is a self-shinker with $|A|^2\equiv 1+\f 56=\f {11}6$.
\end{exam}

It is an interesting question that whether the conditon ``with parallel principal normal'' is necessary in Theorem 1.3-1.6.

For self-shrinkers with the higher codimension, the normal bundle is complicated, which would influence the submanifold properties. Now we consider the simplest case, i.e., the normal bundle is flat. We will prove the following gap theorem:
\begin{thm}\label{thm1-6}
Let $x:M^n\ra\R^{n+p}$ be a complete immersed self-shrinker without boundary and with polynomial volume growth, assume

\vskip 2mm
(i) flat normal bundle, that is,  $R_{\alpha\beta ij}=0$,

(ii) $\sigma_{\alpha\beta}=\f{1}{p}|A|^2\delta_{\alpha\beta}$, where $\sigma_{\alpha\beta}=\sum\limits_{i,j}h^{\alpha}_{ij}h^{\beta}_{ij}$.

\noindent
If the second fundamental form satisfies
\begin{equation}\label{in1-8}
    0\leq |A|^2\leq p,
\end{equation}
then $|A|^2\equiv 0$ and $M^n$ is the hyperplane, or $|A|^2\equiv p$ and $M^n$ is
\begin{equation*}
    x:M^n=N^{mp}\times\R^{n-mp}\hookrightarrow\R^{n+p},
\end{equation*}
where
\begin{equation*}
    N^{mp}=\S^m(\sqrt{m})\times\cdots\times\S^m(\sqrt{m})\hookrightarrow\S^{(m+1)p-1}(\sqrt{mp}).
\end{equation*}

\end{thm}

\begin{rem} As noted in \cite{Smo}, let $\Gamma_1,\Gamma_2,\cdots,\Gamma_m$ be the Abresch-Langer curves, then $\Gamma_1\times\cdots\Gamma_m$ in $\R^{2m}$ is a self-shrinker with
$|{\bf H}|>0$ and flat normal bundle, so the condition \eqref{in1-8} of Theorem \ref{thm1-6} is necessary.
\end{rem}

Finally, we consider the closed self-shrinkers with arbitrary codimension. We have the following simple result.
\begin{prop}\label{prop1-8}
Let $M^n$ be a closed self-shrinker in $\R^{n+p}$, if one of the followings satisfies:
\begin{itemize}
  \item[(1)] $|H|\neq 0$, $\nu=\f H{|H|}$ is parallel in the normal bundle;
  \item[(2)] $|H|^2=const$, or $|H|^2\leq n$, or $|H|^2\geq n$;
  \item[(3)] $|x|^2=const$, or $|x|^2\leq n$, or $|x|^2\geq n$,
\end{itemize}
then $M$ is a minimal submanifold in $\S^{n+p-1}(\sqrt{n})$.
\end{prop}

Note that the condition (1) and (2) imply the closed self-shrinker to be a minimal submanifold in sphere have been proved by Smoczyk \cite{Smo} and Cao-Li \cite{CL2}. In section 8, we will prove that the condition (3) can also imply the self-shrinker is a minimal submanifold in sphere. Then by applying the well-known theorems on the minimal submanifolds in sphere by Ejiri \cite{Ejiri}, H. Li\cite{Li93}, Itoh \cite{Itoh,Itoh2} and Yau \cite{Yau}, Proposition \ref{prop1-8} will imply three simple characterizations for closed self-shrinkers, see Theorem \ref{thm8-1} - \ref{thm8-3}.


\section{Preliminaries}
Let $x:M\to \R^{n+p}$ be an  $n$-dimensional  submanifold of an $(n+p)$-dimensional Euclidean space $R^{n+p}$. Let $\{e_1,\cdots,e_{n}\}$
be a local orthonormal basis of $M$ with respect to the induced
metric, and $\{\theta_1,\cdots, \theta_n\}$ be their dual 1-forms. Let
$e_{n+1},\cdots, e_{n+p}$ be the local unit orthonormal normal
vector fields. In this paper we make the following conventions on the range of
indices:
\begin{equation*}
1\leq i,j,k\leq n; \qquad
n+1\leq \alpha,\beta,\gamma\leq n+p.
\end{equation*}
Then we have the following structure equations (see \cite{CL1,CL2,Li})
\begin{equation*}
dx=\sum\limits_i \theta_i e_i,
\end{equation*}

\begin{equation*}
de_i=\sum\limits_j\theta_{ij}e_j+\sum\limits_{\alpha,j}
h^\alpha_{ij} \theta_je_\alpha,
\end{equation*}

\begin{equation*}
de_\alpha=-\sum\limits_{i,j}h^\alpha_{ij}\theta_j e_i+\sum
\limits_\beta \theta_{\alpha\beta}e_\beta.
\end{equation*}
where  $h^\alpha_{ij}$ denote the the components of the second
fundamental form of $M$. We denote $|A|^2=\sum\limits_{\alpha,i,j}(h^\alpha_{ij})^2$ is the norm square
of the second fundamental form, ${\bf
H}=\sum\limits_{\alpha}H^\alpha e_\alpha=\sum\limits_{\alpha}
(\sum\limits_i h^\alpha_{ii})e_\alpha$ is the mean curvature vector
field, and $H=|\bf{H}|$ is the mean curvature of $M$.

The Gauss equations are given by
\begin{equation}\label{gauss-1}
R_{ijkl}=\sum_{\alpha}(h_{ik}^{\alpha}h_{jl}^{\alpha}-h_{il}^{\alpha}h_{jk}^{\alpha}),
\end{equation}
\begin{equation}
R_{ik}=\sum_{\alpha}H^{\alpha}h_{ik}^{\alpha}-\sum_{\alpha,j}h_{ij}^{\alpha}h_{jk}^{\alpha}.
\end{equation}

The Codazzi equations are given by

\begin{equation}\label{codazzi}
h^\alpha_{ij,k}=h^\alpha_{ik,j},
\end{equation}
where the covariant derivative of $h^\alpha_{ij}$ is defined by
\begin{equation}\label{2-8}
\sum_{k}h^\alpha_{ij,k}\theta_k=dh^\alpha_{ij}+\sum_k
h^\alpha_{kj}\theta_{ki} +\sum_k h^\alpha_{ik}\theta_{kj}+\sum_\beta
h^\beta_{ij}\theta_{\beta\alpha}.
\end{equation}

If we denote by $R_{\alpha\beta i j}$ the  curvature tensor of
the normal connection $\theta_{\alpha\beta}$ in the normal bundle of $x:
M \rightarrow \R^{n+p}$, then the Ricci equations are
\begin{equation}
R_{\alpha\beta i j}=\sum_{k}(h_{ik}^{\alpha}h_{kj}^{\beta}-h_{jk}^{\alpha}h_{ki}^{\beta}).
\end{equation}

By exterior differentiation of $\eqref{2-8}$, we have the following Ricci
identities:
\begin{equation}\label{Ricci}
h^\alpha_{ij,kl}-h^\alpha_{ij,lk}=\sum\limits_m
h^\alpha_{mj}R_{mikl} +\sum\limits_m
h^\alpha_{im}R_{mjkl}+\sum\limits_\beta h^\beta_{ij} R_{\beta\alpha
kl}.
\end{equation}

We define the first and second covariant derivatives, and Laplacian of the mean curvature vector field ${\bf H}=\sum\limits_\alpha H^\alpha e_\alpha$ in the normal bundle $N(M)$ as follows,

\begin{equation}\label{2-11}
\sum\limits_i H^\alpha_{,i}\theta_i=dH^\alpha+\sum\limits_\beta
H^\beta\theta_{\beta\alpha},
\end{equation}

\begin{equation}\label{2-12}
\sum\limits_j H^\alpha_{,ij}\theta_j=dH^\alpha_{,i}+ \sum\limits_j
H^\alpha_{,j}\theta_{ji}+\sum\limits_\beta
H^\beta_{,i}\theta_{\beta\alpha},
\end{equation}

\begin{equation}
\Delta^\perp H^\alpha=\sum\limits_i H^\alpha_{,ii},\qquad
H^\alpha=\sum\limits_k h^\alpha_{kk}.
\end{equation}

Let $f$ be a smooth function on $M$, we define the
covariant derivatives $f_i$, $f_{ij}$, and the Laplacian of $f$ as
follows
\begin{equation*}
df=\sum_i f_i\theta_i,\qquad \sum_j
f_{ij}\theta_j=df_i+\sum_jf_j\theta_{ji},\qquad \Delta f= \sum_i
f_{ii}.
\end{equation*}

\medskip

Now we assume the submanifold $M^n$ satisfies the self-shrinker equation \eqref{1-1}.  The following equations have been derived in \cite{CL2}.

The self-shrinker equation \eqref{1-1} is
equivalent to
\begin{equation}\label{3-1}
H^\alpha=-<x,e_\alpha>, \quad n+1\leq \alpha\leq n+p.
\end{equation}
The first and second covariant derivative of ${\bf H}$ have the following components:

\begin{align}
H^\alpha_{,i}=&\sum\limits_{j}h^\alpha_{ij}<x,e_j>,\label{3-2}\\
H^\alpha_{,ik}=&\sum\limits_{j}h^\alpha_{ij,k}<x,e_j>+h^\alpha_{ik}-\sum\limits_{\beta,j}H^\beta h^\alpha_{ij}h^\beta_{jk},\label{3-3}\\
\Delta H^\alpha=&\sum\limits_{j}H^\alpha_{,j}<x,e_j>+H^\alpha-\sum\limits_{\beta,i,j}H^\beta h^\alpha_{ij}h^\beta_{ij}.\label{LapH}
\end{align}

If ${\bf H}\neq 0$, we can choose local orthogonal frame  $\{e_{\alpha}\}$  for the normal
bundle $NM$ such that $e_{n+p}$ is parallel to the mean curvature vector $\bf H$, that is,

\begin{equation*}
e_{n+p}=\frac{\bf H}{|{\bf H}|}=\nu,\quad H^{n+p}=H=|{\bf H}|,
\qquad H^{\alpha}=0, \quad \alpha\not=n+p.
\end{equation*}

\begin{lem}\label{lem3-1}
Let $x:M^n\ra \R^{n+p}$ be an n-dimensional self-shrinker with $H>0$, $\nabla^{\perp}\nu=0$, then
\begin{equation*}
H^{\alpha}_{,i}=0,\quad H^{\alpha}_{,ij}=0, \alpha\neq n+p,\quad
\text{and} \quad H^{n+p}_{,i}=H_i,\quad H^{n+p}_{,ij}=H_{ij}.
\end{equation*}
\end{lem}
\begin{proof}
Since $\nabla^{\perp}\nu=0$, we have $\theta_{(n+p) \beta}=0$,
from the definition of $H^{\alpha}_{,i}$ in \eqref{2-11}, we have
$H^{\alpha}_{,i}=0$. Then $H^{\alpha}_{,ij}=0$ follows immediately
from \eqref{2-12}.

From ${\bf H}=H^{n+p}e_{n+p}=He_{n+p}$ and
$\nabla^{\perp}e_{n+p}=0$, we have $\nabla_i^{\perp}{\bf
H}=H_ie_{n+p}$, therefore $H^{n+p}_{,i}=H_i$.

From \eqref{2-12}
\begin{eqnarray*}
H^{n+p}_{,ij}\theta_j&=&dH^{n+p}_{,i}+H^{n+p}_{,j}\theta_{ji}+H^{\beta}_{,i}\theta_{\beta
(n+p)}\\
&=&dH_i+H_j\theta_{ji}\\
&=&H_{ij}\theta_j.
\end{eqnarray*}
So we have $H^{n+p}_{,ij}=H_{ij}$.
\end{proof}

\begin{lem}\label{lem3-2}
Let $x:M^n\ra \R^{n+p}$ be an n-dimensional self-shrinker with $H>0$, $\nabla^{\perp}\nu=0$, then
\begin{equation*}
R_{(n+p) \beta
ij}=h^{n+p}_{ik}h^{\beta}_{kj}-h^{\beta}_{ik}h^{n+p}_{kj}=0,\qquad
\sum\limits_{i,j}h^{\alpha}_{ij}h^{n+p}_{ij}=0, \quad\alpha\neq n+p.
\end{equation*}
\end{lem}
\begin{proof}
Since $\nabla^{\perp}\nu=0$ implies $\theta_{(n+p) \beta}=0$, then
$R_{n+p \beta ij}=0$ follows immediately from
\begin{equation*}
d\theta_{(n+p) \beta}-\theta_{(n+p)
\gamma}\wedge\theta_{\gamma\beta}=-\f{1}{2}R_{(n+p)\beta
ij}\theta_i\wedge\theta_j.
\end{equation*}

If $\alpha\neq n+p$, then $H^{\alpha}_{,ij}=0$. From \eqref{LapH}, we have
\begin{equation*}
0=-H^{\alpha}_{,j}<x,e_j>-H^{\alpha}+H^{\beta}h^{\alpha}_{ij}h^{\beta}_{ji}=H\sum\limits_{i,j}h^{\alpha}_{ij}h^{n+p}_{ij}.
\end{equation*}
\end{proof}

\begin{lem}\label{lem3-3}
Let $x:M^n\ra \R^{n+p}$ be an n-dimensional self-shrinker with $H>0$, $\nabla^{\perp}\nu=0$. Denote $|Z|^2=\sum\limits_{i,j}(h^{n+p}_{ij})^2$, we have
\begin{equation*}
\Delta
h^{n+p}_{ij}=\sum\limits_kh^{n+p}_{ij,k}<x,e_k>+h^{n+p}_{ij}-|Z|^2h^{n+p}_{ij}.
\end{equation*}
\end{lem}
\begin{proof} By use of \eqref{codazzi}, \eqref{Ricci}, Lemma \ref{lem3-2}, \eqref{gauss-1} and \eqref{3-3}, we have the following calculations:
\begin{eqnarray*}
\Delta h^{n+p}_{ij}&=&h^{n+p}_{ij,kk}=h^{n+p}_{ik,jk}\\
&=&h^{n+p}_{kk,ij}+h^{n+p}_{mk}R_{mijk}+h^{n+p}_{mi}R_{mkjk}+h^{\beta}_{ik}R_{\beta(n+p)
jk}\\
&=&H^{n+p}_{,ij}+h^{n+p}_{mk}(h^{\alpha}_{mj}h^{\alpha}_{ik}-h^{\alpha}_{mk}h^{\alpha}_{ij})+h^{n+p}_{mi}(h^{\alpha}_{mj}H^{\alpha}-h^{\alpha}_{mk}h^{\alpha}_{kj})\\
&=&H^{n+p}_{,ij}-h^{n+p}_{mk}h^{\alpha}_{mk}h^{\alpha}_{ij}+h^{n+p}_{mi}h^{n+p}_{mj}H^{n+p}\\
&=&h^{n+p}_{ij,k}<x,e_k>+h^{n+p}_{ij}-h^{n+p}_{mk}h^{n+p}_{mk}h^{n+p}_{ij}\\
&=&h^{n+p}_{ij,k}<x,e_k>+h^{n+p}_{ij}-|Z|^2h^{n+p}_{ij}.
\end{eqnarray*}
\end{proof}

Concerning the term $|Z|^2$, we have the following inequalities which holds for all submanifolds with parallel principal normal (may not be a self-shrinker).
\begin{lem}\label{lem3-5}
If we fix a point $q$ and choose a frame $e_i$, $i=1,\ldots,n$, such
that $h^{n+p}_{ij}$ is diagonal at $q$, i.e.
$h^{n+p}_{ij}=\lambda_i\delta_{ij}$,  then we have at $q$ that
\begin{align}
&|\nabla |Z||^2\leq \sum_{i,k}(h^{n+p}_{ii,k})^2\leq\sum_{i,j,k}(h^{n+p}_{ij,k})^2\label{3-12}\\
&(1+\f{2}{n+1})|\nabla|Z||^2\leq\sum_{i,j,k}(h^{n+p}_{ij,k})^2+\f{2n}{n+1}|\nabla H|^2.\label{3-13}
\end{align}
\end{lem}
\begin{proof}
Since $h^{n+p}_{ij}$ is symmetric, we can choose $e_i$,
$i=1,\ldots,n$, such that $h^{n+p}_{ij}=\lambda_i\delta_{ij}$ at the
fixed point $q$. By $\nabla |Z|^2=2|Z|\nabla|Z|$, we have at $q$
\begin{equation*}
4|Z|^2|\nabla|Z||^2=\sum_k(2h^{n+p}_{ij}h^{n+p}_{ij,k})^2=4\sum_k(\sum_i\lambda_ih^{n+p}_{ii,k})^2\leq
4|Z|^2\sum_{i,k}(h^{n+p}_{ii,k})^2,
\end{equation*}
where the inequality used the Cauchy-Schwarz inequality, this proves
\eqref{3-12}.

To show \eqref{3-13}, we have by \eqref{3-12}
\begin{eqnarray*}
|\nabla|Z||^2&\leq&\sum_{i,k}(h^{n+p}_{ii,k})^2\\
&=&\sum_{i\neq k}(h^{n+p}_{ii,k})^2+\sum_i(h^{n+p}_{ii,i})^2\\
&=&\sum_{i\neq k}(h^{n+p}_{ii,k})^2+\sum_i(H^{n+p}_{,i}-\sum_{j\neq
i}h^{n+p}_{jj,i})^2\\
&\leq&\sum_{i\neq
k}(h^{n+p}_{ii,k})^2+n\sum_i((H^{n+p}_{,i})^2+\sum_{j\neq
i}(h^{n+p}_{jj,i})^2)\\
&=&n|\nabla H|^2+(n+1)\sum_{i\neq k}(h^{n+p}_{ii,k})^2\\
&=&n|\nabla H|^2+\f{n+1}{2}(\sum_{i\neq
k}(h^{n+p}_{ik,i})^2+\sum_{i\neq k}(h^{n+p}_{ki,i})^2),
\end{eqnarray*}
where we used the algebraic fact $(\sum\limits_{i=1}^na_i)^2\leq n\sum\limits_{i=1}^na_i^2$ in the second inequality, and  $H^{n+p}_{,i}=H_i$,
Codazzi equations in the last two equalities. Thus we have
\begin{eqnarray*}
(1+\f{2}{n+1})|\nabla|Z||^2&\leq&\f{2n}{n+1}|\nabla
H|^2+\sum_{i,k}(h^{n+p}_{ii,k})^2+\sum_{i\neq
k}(h^{n+p}_{ik,i})^2+\sum_{i\neq k}(h^{n+p}_{ki,i})^2\\
&\leq&\f{2n}{n+1}|\nabla H|^2+\sum_{i,j,k}(h^{n+p}_{ij,k})^2,
\end{eqnarray*}
which completes the proof.
\end{proof}

\section{Some integral estimates}

Recall the following operator $\cal{L}$ which was
introduced and studied firstly on self-shrinkers by Colding and
Minicozzi (see (3.7) in \cite{CM2}):
\begin{equation*}
\cal{L}=\Delta-<x,\nabla(\cdot)>=e^{\f{|x|^2}{2}}\textrm{div}(e^{-\f{|x|^2}{2}}\nabla\cdot),
\end{equation*}
where $\Delta$, $\nabla$ and $\textrm{div}$ denote the Laplacian, gradient and divergent operator on the self-shrinker respectively, $<\cdot,\cdot>$ denotes the standard inner product in $\R^{n+p}$. The operator $\cal{L}$ is self-adjoint in a weighted $L^2$ space. The next two results were proved by Colding-Minicozzi \cite{CM2} for hypersurface self-shrinkers but can be stated in the same way for self-shrinkers in arbitrary codimension.
\begin{lem}\label{lem2-1}
If $x: M^n\ra \R^{n+p}$ is a submanifold, $u$ is a $C^1$ function
with compact support, and $v$ is a $C^2$ function, then
\begin{equation}
\int_Mu(\cal{L}v)e^{-\f{|x|^2}{2}}=-\int_M<\nabla v,\nabla
u>e^{-\f{|x|^2}{2}}.
\end{equation}
\end{lem}
\begin{corr}\label{cor2-1}
Suppose that $x: M^n\ra\R^{n+p}$ is a complete submanifold without boundary, if $u, v$ are $C^2$ functions satisfying
\begin{equation*}
\int_M(|u\nabla v|+|\nabla u||\nabla
v|+|u\cal{L}v|)e^{-\f{|x|^2}{2}}<+\infty,
\end{equation*}
then we get
\begin{equation}
\int_Mu(\cal{L}v)e^{-\f{|x|^2}{2}}=-\int_M<\nabla v,\nabla
u>e^{-\f{|x|^2}{2}}.
\end{equation}
\end{corr}

Using the operator $\mathcal{L}$, Lemma 2.3 has the following Corollary,
\begin{corr}\label{cor3-1}
Let $x:M^n\ra \R^{n+p}$ be an n-dimensional complete self-shrinker with $H>0$, $\nabla^{\perp}\nu=0$, then
\begin{eqnarray*}
\f{1}{2}\cal{L}
|Z|^2&=&|Z|^2-|Z|^4+\sum_{i,j,k}(h^{n+p}_{ij,k})^2,\\
\cal{L}|Z|&=&|Z|-|Z|^3+\f{\sum_{i,j,k}(h^{n+p}_{ij,k})^2}{|Z|}-\f{|\nabla
|Z||^2}{|Z|}.
\end{eqnarray*}
\end{corr}
\begin{rem}
We note that our assumption ``$H>0$" implies ``$|Z|>0$" because of $|Z|^2\geq H^2/n$.
\end{rem}

\begin{lem}\label{lem3-4}
Let $x:M^n\ra \R^{n+p}$ be an n-dimensional complete self-shrinker with $H>0$, $\nabla^{\perp}\nu=0$, then
\begin{eqnarray*}
\cal{L}H&=&H-|Z|^2H,\\
\cal{L}\log{H}&=&1-|Z|^2-|\nabla\log{H}|^2.
\end{eqnarray*}
\end{lem}
\begin{proof}
The two equations just follow from \eqref{LapH} and Lemma \ref{lem3-1}.
\end{proof}

\begin{lem}\label{lem3-6}
Let $x:M^n\ra \R^{n+p}$ be an n-dimensional complete self-shrinker with $H>0$, $\nabla^{\perp}\nu=0$. If $\varphi$ is in the weighted $W^{1,2}$ space, i.e.
\begin{equation*}
\int_M(|\varphi|^2+|\nabla \varphi|^2)e^{-\f{|x|^2}{2}}<+\infty,
\end{equation*}
then
\begin{equation}
\int_M\varphi^2(|Z|^2+\f{1}{2}|\nabla\log H|^2)e^{-\f{|x|^2}{2}}\leq
\int_M(2|\nabla\varphi|^2+\varphi^2)e^{-\f{|x|^2}{2}}.
\end{equation}
\end{lem}
\begin{proof}
Suppose that $\eta$ is a function with compact support, from the
self-adjointness of $\cal{L}$  and Lemma 3.4 we have
\begin{eqnarray*}
\int_M<\nabla\eta^2,\nabla\log
H>e^{-\f{|x|^2}{2}}&=&-\int_M\eta^2(\cal{L}\log
H)e^{-\f{|x|^2}{2}}\\
&=&\int_M\eta^2(|Z|^2-1+|\nabla\log H|^2)e^{-\f{|x|^2}{2}}.
\end{eqnarray*}
Combining this with the Cauchy-Schwarz inequality
\begin{eqnarray*}
<\nabla\eta^2,\nabla\log H>&\leq&
2|\nabla\eta|^2+\f{1}{2}\eta^2|\nabla\log H|^2
\end{eqnarray*}
gives that
\begin{equation*}
\int_M\eta^2(|Z|^2+\f{1}{2}|\nabla\log H|^2)e^{-\f{|x|^2}{2}}\leq
\int_M(2|\nabla\eta|^2+\eta^2)e^{-\f{|x|^2}{2}}.
\end{equation*}
Now we choose a sequence of cut-off function $\eta_j$ which satisfies
\begin{equation*}
    \eta_j=\left\{\begin{array}{ll}
                    1, & \textrm{in }B_j \\
                    0, & \textrm{outside }B_{j+1}
                  \end{array}\right.,\quad 0\leq \eta_j\leq 1, \quad |\nabla\eta_j|\leq C,
\end{equation*}
where $B_j=M\cap B_j(0)$ with $B_j(0)$ is the Euclidean ball of radius $j$ centered at the origin. Applying the above inequality
with $\eta=\eta_j\varphi$, letting $j\ra \infty$, and using the
dominated convergence theorem, we complete the proof of the Lemma.
\end{proof}

The next proposition gives weighted estimates for the principle normal second fundamental form and its covariant derivatives.
\begin{prop}\label{lem3-7}
Let $x:M^n\ra \R^{n+p}$ be an n-dimensional complete self-shrinker with $H>0$, $\nabla^{\perp}\nu=0$. If $M^n$ has polynomial volume growth, then
\begin{equation}
\int_M(|Z|^2+|Z|^4+|\nabla|Z||^2+\sum_{ijk}(h^{n+p}_{ij,k})^2)e^{-\f{|x|^2}{2}}<\infty.
\end{equation}
\end{prop}
\begin{proof}
For any compactly supported function $\varphi$, self-adjointness
of $\cal{L}$ and Lemma 3.4 imply
\begin{eqnarray*}
\int_M<\nabla\varphi^2,\nabla\log
H>e^{-\f{|x|^2}{2}}&=&-\int_M\varphi^2(\cal{L}\log
H)e^{-\f{|x|^2}{2}}\\
&=&\int_M\varphi^2(|Z|^2-1+|\nabla\log H|^2)e^{-\f{|x|^2}{2}}.
\end{eqnarray*}
Combining this with the Cauchy-Schwarz inequality
\begin{equation*}
<\nabla \varphi^2, \nabla\log H>~\leq
~|\nabla\varphi|^2+\varphi^2|\nabla\log H|^2
\end{equation*}
gives the following stability inequality (cf. \cite{SSY})
\begin{equation*}
\int_M\varphi^2|Z|^2e^{-\f{|x|^2}{2}}~\leq~\int_M(|\nabla\varphi|^2+\varphi^2)e^{-\f{|x|^2}{2}}.
\end{equation*}
Let $\varphi=\eta|Z|$, where $\eta\geq 0$ has compact support, for $\epsilon>0$, we have
\begin{equation}\label{3-30}
\begin{split}
\int_M\eta^2|Z|^4e^{-\f{|x|^2}{2}}~\leq&
\int_M\left(\eta^2|\nabla|Z||^2+2\eta|Z||\nabla\eta||\nabla|Z||\right.\\
&\qquad+\left.|Z|^2|\nabla\eta|^2+\eta^2|Z|^2\right)e^{-\f{|x|^2}{2}}\\
\leq&
(1+\epsilon)\int_M\eta^2|\nabla|Z||^2e^{-\f{|x|^2}{2}}\\
&\qquad+\int_M|Z|^2((1+\f{1}{\epsilon})|\nabla\eta|^2+\eta^2)e^{-\f{|x|^2}{2}}.
\end{split}
\end{equation}
Corollary 3.3 and Lemma 2.4 give the
inequality
\begin{eqnarray*}
\cal{L}|Z|^2&\geq&2(1+\f{2}{n+1})|\nabla|Z||^2-\f{4n}{n+1}|\nabla
H|^2+2|Z|^2-2|Z|^4.
\end{eqnarray*}
Integrating this with $\f{1}{2}\eta^2$, it follows from the
self-adjointness of $\cal{L}$ that
\begin{align*}
-2\int_M\eta|Z|<\nabla\eta,\nabla|Z|>e^{-\f{|x|^2}{2}}\geq&\int_M\left(\eta^2(1+\f{2}{n+1})|\nabla|Z||^2\right.\\
&\qquad \left.-\f{2n}{n+1}\eta^2|\nabla H|^2-\eta^2|Z|^4\right)e^{-\f{|x|^2}{2}}.
\end{align*}
Using the inequality $2ab\leq \epsilon a^2+\f{b^2}{\epsilon}$ gives
\begin{eqnarray}
&&\int_M(\eta^2|Z|^4+\f{2n}{n+1}\eta^2|\nabla
H|^2+\f{1}{\epsilon}|Z|^2|\nabla\eta|^2)e^{-\f{|x|^2}{2}}\label{3-33}\\
&\geq&\int_M(1+\f{2}{n+1}-\epsilon)\eta^2|\nabla|Z||^2e^{-\f{|x|^2}{2}}.\nonumber
\end{eqnarray}
Assume $|\eta|\leq 1$ and $|\nabla\eta|\leq 1$, combining
\eqref{3-30} and \eqref{3-33} gives
\begin{align*}
\int_M\eta^2|Z|^4e^{-\f{|x|^2}{2}}~\leq~&\f{1+\epsilon}{1+\f{2}{n+1}-\epsilon}\int_M\eta^2|Z|^4e^{-\f{|x|^2}{2}}\\
&\qquad+C_{\epsilon}\int_M(|\nabla
H|^2+|Z|^2)e^{-\f{|x|^2}{2}}.
\end{align*}
Choose $\epsilon>0$ small, such that
$\f{1+\epsilon}{1+\f{2}{n+1}-\epsilon}<1$, then we have
\begin{eqnarray}
\int_M\eta^2|Z|^4e^{-\f{|x|^2}{2}}&\leq& C\int_M(|\nabla
H|^2+|Z|^2)e^{-\f{|x|^2}{2}}\label{3-35}\\
&\leq& C\int_M|Z|^2(1+|x|^2)e^{-\f{|x|^2}{2}},\nonumber
\end{eqnarray}
where the second inequality is due to \eqref{3-2} and Lemma 2.1. Since $H>0$, Lemma 3.5 and the polynomial
volume growth give that
$\int_M|Z|^2(1+|x|^2)e^{-\f{|x|^2}{2}}<\infty$, thus \eqref{3-35}
and the dominated convergence theorem give that
$\int_M|Z|^4e^{-\f{|x|^2}{2}}<\infty$, then
$\int_M|\nabla|Z||^2e^{-\f{|x|^2}{2}}<\infty$ follows immediately
from \eqref{3-33} and the dominated convergence theorem.

\smallskip

To show
$\int_M\sum\limits_{i,j,k}(h^{n+p}_{ij,k})^2e^{-\f{|x|^2}{2}}<\infty$,
we integrate the first equation in Corollary 3.3 with
$\eta^2$, the self-adjointness of $\cal{L}$ implies
\begin{eqnarray*}
\int_M\eta^2\sum_{ijk}(h^{n+p}_{ij,k})^2e^{-\f{|x|^2}{2}}&=&\int_M\eta^2(|Z|^4-|Z|^2)e^{-\f{|x|^2}{2}}\\
&&\quad-\int_M2\eta|Z|<\nabla\eta,\nabla|Z|>e^{-\f{|x|^2}{2}}\\
&\leq&\int_M(\eta^2|Z|^4+|\nabla\eta|^2|\nabla|Z||^2)e^{-\f{|x|^2}{2}}~<\infty.
\end{eqnarray*}
The dominated convergence theorem gives that
$\int_M\sum\limits_{ijk}(h^{n+p}_{ij,k})^2e^{-\f{|x|^2}{2}}<\infty$,
this completes the proof.
\end{proof}

Now we will prove the weighted integral estimates that will be needed in the next section, to guarantee the self-adjointness of $\mathcal{L}$ can apply on complete self-shrinkers.
\begin{prop}\label{prop3-7}
Let $x:M^n\ra \R^{n+p}$ be an n-dimensional complete self-shrinker with $H>0$, $\nabla^{\perp}\nu=0$. If $M^n$ has polynomial volume growth, then
\begin{align}
&\int_M\left(|Z|^2|\nabla\log H|+|\nabla|Z|^2||\nabla\log
H|+|Z|^2|\cal{L}\log H|\right)e^{-\f{|x|^2}{2}}<\infty,\label{3-22}\\
&\int_M(|Z||\nabla|Z||+|\nabla|Z||^2+|Z||\cal{L}|Z||)e^{-\f{|x|^2}{2}}<\infty.\label{3-23}
\end{align}
\end{prop}

\begin{proof}
Proposition 3.6 implies that $|Z|$ is in the weighted $W^{1,2}$ space, so
Lemma 3.5 gives that
\begin{equation*}
\int_M|Z|^2|\nabla\log H|^2e^{-\f{|x|^2}{2}}~<\infty.
\end{equation*}
Then
\begin{equation*}
\int_M|\nabla|Z|^2||\nabla\log H|e^{-\f{|x|^2}{2}}~\leq
\int_M(|\nabla|Z||^2+|Z|^2|\nabla\log H|^2)e^{-\f{|x|^2}{2}}<\infty.
\end{equation*}
From Lemma 3.4, we have
\begin{equation*}
\int_M|Z|^2|\cal{L}\log
H|e^{-\f{|x|^2}{2}}=\int_M|Z|^2|1-|Z|^2-|\nabla\log
H|^2|e^{-\f{|x|^2}{2}}<\infty.
\end{equation*}
This gives the first part of Proposition 3.7; From Corollary 3.3
\begin{equation*}
\int_M|Z||\cal{L}|Z||e^{-\f{|x|^2}{2}}=\int_M(|Z|^2-|Z|^4+\sum_{ijk}(h^{n+p}_{ij,k})^2-|\nabla|Z||^2)e^{-\f{|x|^2}{2}}<\infty.
\end{equation*}
So the second part follows from Proposition 3.6.
\end{proof}

\section{Proof of Theorem 1.1}

In this section, we will give the proof of Theorem 1.1. First we prove two geometric identities, which is the key for proving the classification.
\begin{lem}
Let $x:M^n\ra \R^{n+p}$ be an n-dimensional complete self-shrinker with $H>0$, $\nabla^{\perp}\nu=0$. If $M^n$ has polynomial volume growth, then
\begin{eqnarray}
|Z|&=&\beta H\quad \textrm{for some positive constant } \beta,\label{4-1}\\
|\nabla|Z||^2&=&\sum_{i,j,k}(h^{n+p}_{ij,k})^2.\label{4-2}
\end{eqnarray}
\end{lem}
\begin{proof}
By \eqref{3-22}, we can apply Corollary 3.2 to $|Z|^2$ and $\log
H$ to get by use of Lemma 3.4
\begin{equation}\label{4-8}
\begin{split}
\int_M<\nabla|Z|^2,\nabla\log
H>e^{-\f{|x|^2}{2}}~&=-\int_M|Z|^2(\cal{L}\log H)e^{-\f{|x|^2}{2}}\\
&=\int_M|Z|^2(|Z|^2-1+|\nabla\log H|^2)e^{-\f{|x|^2}{2}}
\end{split}
\end{equation}
Similarly, by \eqref{3-23}, we apply Corollary 3.2 to two copies
of $|Z|$ to get by use of Corollary 3.3 and (2.14)
\begin{equation}\label{4-9}
\begin{split}
\int_M|\nabla|Z||^2e^{-\f{|x|^2}{2}}~&=-\int_M|Z|\cal{L}|Z|e^{-\f{|x|^2}{2}}\\
&=\int_M\left(|Z|^4-|Z|^2+|\nabla|Z||^2-\sum_{i,j,k}(h^{n+p}_{ij,k})^2\right)e^{-\f{|x|^2}{2}}\\
&\leq \int_M\left(|Z|^4-|Z|^2\right)e^{-\f{|x|^2}{2}}.
\end{split}
\end{equation}
Combining
\eqref{4-8} and \eqref{4-9} give
\begin{align*}
0&\geq \int_M\left(|\nabla|Z||^2-2|Z|<\nabla|Z|,\nabla\log
H>+|Z|^2|\nabla\log H|^2\right)e^{-\f{|x|^2}{2}}\\
&=\int_M\left|\nabla|Z|-|Z|\nabla\log H\right|^2e^{-\f{|x|^2}{2}}
\end{align*}
So we conclude that $\nabla|Z|\equiv|Z|\nabla\log H$, therefore,
$|Z|=\beta H$ for some constant $\beta>0$. And the inequality in
\eqref{4-9} must be equality,  so we have
$|\nabla|Z||^2=\sum\limits_{ijk}(h^{n+p}_{ij,k})^2$.
\end{proof}

\begin{lem}\label{two-cases}
Let $x:M^n\ra \R^{n+p}$ be an n-dimensional complete connected self-shrinker with $H>0$, $\nabla^{\perp}\nu=0$. If $M^n$ has polynomial volume growth, then one of the following two cases holds

\medskip

(i) $\nabla^\perp{\bf H}\equiv 0$ and $h^{n+p}_{ij,k}\equiv 0$,

\vskip 1mm

(ii) $(h^{n+p}_{ij})$ admits only one nonzero eigenvalue $H$, in this case, $|Z|^2=H^2$.

\end{lem}

\begin{proof}
As in Lemma 2.4, we fix a point $q$ and choose a frame $e_i$,
$i=1,\ldots,n$, such that $h^{n+p}_{ij}$ is diagonal at $q$, i.e.
$h^{n+p}_{ij}=\lambda_i\delta_{ij}$ then we have at $q$ that
\begin{equation*}
|Z|^2|\nabla |Z||^2=\sum_k(\sum_ih^{n+p}_{ii,k}\lambda_i)^2\leq
|Z|^2\sum_{i,k}(h^{n+p}_{ii,k})^2\leq
|Z|^2\sum_{i,j,k}(h^{n+p}_{ij,k})^2.
\end{equation*}
By \eqref{4-2}, the above two inequalities must be equalities, so we have:
\begin{enumerate}
\item[(i)] For each $k$, there exists a constant $\alpha_k$ such
that $h^{n+p}_{ii,k}=\alpha_k\lambda_i$ for every $i$.
\item[(ii)] If $i\neq j$, then $h^{n+p}_{ij,k}=0$.
\end{enumerate}
By the Codazzi equation, $(ii)$ implies that
\begin{enumerate}
\item[(ii)'] $h^{n+p}_{ij,k}=0$ unless $i=j=k$.
\end{enumerate}
If $\lambda_i\neq 0$, for $j\neq i$, $0=h^{n+p}_{ii,j}=\alpha_j\lambda_i$ so we have $\alpha_j=0$. If the rank of $(h^{n+p}_{ij})$ is at least two at $q$, then $\alpha_j=0$ for all $j\in\{1,\cdots,n\}$. Thus (i),(ii) imply $H^{n+p}_{,k}=0$ for all $k\in\{1,\cdots,n\}$ and $h^{n+p}_{ij,k}=0$ for all $i,j,k\in \{1,\cdots,n\}$. If the rank of $(h^{n+p}_{ij})$ is one at $q$, then $H$ is the only nonzero eigenvalue of $(h^{n+p}_{ij})$, and then $|Z|^2=H^2$.

Next we will show that if the rank of $(h^{n+p}_{ij})$ is at least two at some $q$, then the rank of $(h^{n+p}_{ij})$ is at least two everywhere. For each $x\in M$, let $\lambda_1(x)$ and $\lambda_2(x)$ be the two eigenvalues of $(h^{n+p}_{ij}(x))$ that are largest in absolute value and define the set
\begin{equation*}
    \Omega=\{x\in M|\lambda_1(x)=\lambda_1(q),\quad \lambda_2(x)=\lambda_2(q)\}.
\end{equation*}
Then $\Omega$ is nonempty since $q\in \Omega$. $\lambda_i(x)$ are continuous in $x$, so $\Omega$ is closed. For any $x\in\Omega$, the rank of $(h^{n+p}_{ij})$ is at least two at $x$, this is an open condition, so there is an open set $U$ containing $x$ where the rank of $(h^{n+p}_{ij})$ is at least two, then $h^{n+p}_{ij,k}\equiv 0$ on this set $U$, and the eigenvalues of $(h^{n+p}_{ij})$ are constant on $U$. This implies $U\subset \Omega$, and therefore $\Omega$ is open. Since $M$ is connected, we conclude that $\Omega=M$, therefore the rank of $(h^{n+p}_{ij})$ is at least two everywhere. This implies the Case (i).

Since $H>0$, the remaining case is where the rank of $(h^{n+p}_{ij})$ is exactly one at every point, this implies Case (ii). This completes the proof of Lemma 4.2.
\end{proof}

\vskip 3mm
\begin{proof}[Proof of Theorem 1.1]
We can treat the two cases in Lemma \ref{two-cases} separately by following the argument of K. Smoczyk in \cite{Smo} to complete the proof. Note that in Case I, $h^{n+p}_{ij,k}\equiv 0$ implies $|A^{\nu}|^2=|Z|^2\equiv \textrm{constant}$,  combining with the assumption $|A|^2-|A^\nu|^2\leq c$, we can complete the proof of Case I. While for the Case II, the assumption  $|A|^2-|A^\nu|^2\leq c$ is sufficient for us to complete the proof by following the  Smoczyk's arguments in \cite{Smo}.
\end{proof}

\section{Rigidity of Self-shrinkers in higher codimension}

\subsection{Self-shrinkers with codimension two}
Assume $x: M^n\ra\R^{n+2}$ is a complete imbedded self-shrinker without boundary and with polynomial volume growth in codimension two, we will give the proof of Theorem 1.3 and 1.4.

\vskip 3mm
\begin{proof}[Proof of Theorem 1.3]
Under the condition of the Theorem 1.3, Theorem A and Theorem 1.1 imply
\begin{equation}\label{codim2-1}
    M^n=\Gamma\times \R^{n-1},
\end{equation}
where $\Gamma$ is one of the Abresch-Langer curves, or
\begin{equation}\label{codim2-2}
    M^n=\td{M}^r\times\R^{n-r}\hookrightarrow\S^{r+1}(\sqrt{r})\times\R^{n-r}\hookrightarrow\R^{n+2},\quad 0<r\leq n,
\end{equation}
where $\td{M}^r\hookrightarrow\S^{r+1}(\sqrt{r})$ is a minimal hypersurface. Note that the only simple closed one of the Abresch-Langer curves is the circle, so the first case \eqref{codim2-1} is $M^n=\S^1(1)\times\R^{n-1}$ with $|A|^2\equiv 1$. Then we consider the second case \eqref{codim2-2}. Denote $A$ the second fundamental form of $M^n\ra\R^{n+2}$, and $\td{A}$ the second fundamental form of $\td{M}^r\hookrightarrow\S^{r+1}(\sqrt{r})$. Then
\begin{align}
|A|^2 =& |\td{A}|^2+1,\label{ssf}\\
\f 12\tilde{\Delta}|\td{A}|^2=& |\tilde{\nabla} \td{A}|^2+|\td{A}|^2(1-|\td{A}|^2),\label{eq5-2}
\end{align}
where $\tilde{\Delta}, \tilde{\nabla}$ denote Laplacian and covariant derivatives with respect to the induced metric on $\tilde{M}$. The Simons' equality \eqref{eq5-2} can be found in Simons' paper \cite{Si}.  For convenience of readers, we give a proof here: Denote $\tilde{h}_{ij}, \tilde{R}_{ijkl}$ the components of the second fundamental form and the curvature tensor of $\tilde{M}^r$ in $\S^{r+1}(\sqrt{r})$.  Since the sphere $\S^{r+1}(\sqrt{r})$ has constant curvature $\frac 1r$, the Ricci identities and Gauss-Codazzi equations give that
\begin{align*}
\tilde{\Delta}\tilde{h}_{ij}=&\tilde{h}_{ij,kk}=\tilde{h}_{ki,jk}\\
=&\tilde{h}_{ki,kj}+\tilde{R}_{kjkl}\tilde{h}_{li}+\tilde{R}_{kjil}\tilde{h}_{kl}\\
=&((r-1)\frac 1r\delta_{jl}-\td{h}_{jk}\td{h}_{kl})\td{h}_{li}+(\frac 1r(\delta_{ik}\delta_{jl}-\delta_{kl}\delta_{ji})+\td{h}_{ki}\td{h}_{jl}-\td{h}_{kl}\td{h}_{ji})\td{h}_{kl}\\
=&(1-|\tilde{A}|^2)\tilde{h}_{ij},
\end{align*}
where we used that $\tilde{M}^r\hookrightarrow\S^{r+1}(\sqrt{r})$ is a minimal hypersurface in the last two equalities. Then it follows that
\begin{align*}
\frac 12\tilde{\Delta}|\tilde{A}|^2=\tilde{h}_{ij}\tilde{\Delta}\tilde{h}_{ij}+|\tilde{\nabla}\tilde{A}|^2= |\tilde{\nabla} \td{A}|^2+|\td{A}|^2(1-|\td{A}|^2).
\end{align*}

Since $x:M^n\ra\R^{n+2}$ is a embedded self-shrinker without boundary and with polynomial volume growth, by Cheng-Zhou (see Theorem 4.1 in \cite{ChZ}), $x:M^n
\ra\R^{n+2}$ is proper, thus we have that $\td{M}^r$ is closed. Therefore \eqref{eq5-2} implies that if
\begin{equation}\label{gap-1}
    0\leq|\td{A}|^2\leq 1,
\end{equation}
we have either $|\td{A}|\equiv 0$ and $\td{M}^r$ is totally geodesic in $\S^{r+1}(\sqrt{r})$, that is $\td{M}^r=\S^r(\sqrt{r})$; or $|\td{A}|\equiv 1$ and $\td{M}^r$ is the Clifford minimal hypersurface in $\S^{r+1}(\sqrt{r})$, that is $\td{M}^r=\S^k(\sqrt{k})\times\S^{r-k}(\sqrt{r-k})$. From \eqref{ssf}, the condition \eqref{gap-1} is equivalent to
\begin{equation}\label{gap-2}
   1\leq|A|^2\leq 2.
\end{equation}
So we conclude that if we have \eqref{gap-2}, then there are two possibilities:

\begin{enumerate}
  \item[(1)] $|A|^2\equiv 1$ and $M^n=\S^m(\sqrt{m})\times\R^{n-m}$ with $1\leq m\leq n$;
  \item[(2)] $|A|^2\equiv 2$ and $M^n=\S^k(\sqrt{k})\times\S^{r-k}(\sqrt{r-k})\times\R^{n-r}$ with $2\leq r\leq n$ and $1\leq k\leq r-1$.
\end{enumerate}

\end{proof}

\vskip 3mm
\begin{proof}[Proof of Theorem 1.4]
As the proof of Theorem 1.3, Theorem A and Theorem 1.1 imply $M^n$ must be one of the two cases \eqref{codim2-1} and \eqref{codim2-2}. Note that the case \eqref{codim2-1} has $|A|^2\equiv 1$, which violates with the assumption \eqref{th1-3}, so $M$ must be the case \eqref{codim2-2}, that is
\begin{equation*}
    M^n=\td{M}^r\times\R^{n-r}\hookrightarrow\S^{r+1}(\sqrt{r})\times\R^{n-r}\hookrightarrow\R^{n+2},\quad 0<r\leq n,
\end{equation*}
where $\td{M}^r\hookrightarrow\S^{r+1}(\sqrt{r})$ is a closed minimal hypersurface, then \eqref{ssf} holds. Q. Ding and Y. L. Xin \cite{DX1} proved that there exists a constant $\delta>0$ such that if
\begin{equation}
    1\leq|\td{A}|^2\leq 1+\delta,
\end{equation}
then $|\td{A}|\equiv 1$ and $M$ is the Clifford minimal hypersurface. From \eqref{ssf} again, we conclude that if
\begin{equation}
    2\leq|A|^2\leq 2+\delta,
\end{equation}
we have $|A|^2\equiv 2$ and $M^n=\S^k(\sqrt{k})\times\S^{r-k}(\sqrt{r-k})\times\R^{n-r}$ with $2\leq r\leq n$ and $1\leq k\leq r-1$.
\end{proof}

\subsection{Self-shrinkers of two dimension}

In this subsection, we assume $x: M^2\ra\R^{2+p}$ is a 2-dimensional complete imbedded self-shrinker without boundary and with polynomial volume growth, and we will give the proof of Theorem 1.5 and 1.6.

\begin{proof}[Proof of Theorem 1.5]
Theorem A and Theorem 1.1 imply that
\begin{equation*}
    M^2=\td{M}^1\times\R,\quad \textrm{or}\quad M^2=\td{M}^2,
\end{equation*}
where $\td{M}^1\hookrightarrow\S^{p}(1)$ and $\td{M}^2\hookrightarrow\S^{p+1}(\sqrt{2})$ are minimal submanifolds. Then $\td{M}^1=\S^1(1)$ and $M^2=\td{M}^1\times\R\hookrightarrow\R^{2+p}$ has $|A|^2\equiv 1$. Since $x:M^2\ra\R^{2+p}$ is an embedded self-shrinker without boundary and with polynomial volume growth, $\td{M}^2\hookrightarrow\S^{p+1}(\sqrt{2})$ are closed minimal submanifolds. Denote $\td{A}$ the second fundamental form of $\td{M}^2\hookrightarrow\S^{p+1}(\sqrt{2})$, then
\begin{equation*}
    |A|^2=|\td{A}|^2+1,
\end{equation*}
where $A$ is the second fundamental form of $x:M^2\ra\R^{2+p}$. Then the condition \eqref{dim2} is equivalent to
\begin{equation}
    0\leq |\td{A}|^2\leq \f 23.
\end{equation}
A well-known theorem (cf. Theorem B in \cite{KS}) implies either

\begin{enumerate}
  \item[(1)] $|\td{A}|^2\equiv 0$ and $\td{M}^2=\S^2(\sqrt{2})$. or
  \item[(2)] $|\td{A}|^2\equiv \f 23$ and $\td{M}^2=\S^2(\sqrt{6})\hookrightarrow\S^{4}(\sqrt{2})$ is a Veronese surface.
\end{enumerate}

\noindent Therefore in terms of $|A|^2$, we have two possibilities:

\begin{enumerate}
  \item[(1)] $|A|^2\equiv 1$ and $M^2=\S^2(\sqrt{2})$ in $\R^{3}$
  \item[(2)] $|A|^2\equiv \f 53$ and the self-shrinker has the form $x:M^2=\S^2(\sqrt{6})\ra\S^{4}(\sqrt{2})\hookrightarrow\R^{5}$.
\end{enumerate}

Note that $M^2=\S^1(1)\times\R\hookrightarrow\R^{2+p}$ also has $|A|^2\equiv 1$, this completes the proof of Theorem 1.5.
\end{proof}

\vskip 3mm

\begin{proof}[Proof of Theorem 1.6]
Theorem A and Theorem 1.1 imply that
\begin{equation*}
    M^2=\td{M}^1\times\R,\quad \textrm{or}\quad M^2=\td{M}^2,
\end{equation*}
where $\td{M}^1\hookrightarrow\S^{p}(1)$ and $\td{M}^2\hookrightarrow\S^{p+1}(\sqrt{2})$ are minimal submanifolds. Then $\td{M}^1=\S^1(1)$ and $M^2=\td{M}^1\times\R\hookrightarrow\R^{2+p}$ has $|A|^2\equiv 1$, this is impossible. So
 \begin{equation*}
    M^2=\td{M}^2\hookrightarrow\S^{p+1}(\sqrt{2})
 \end{equation*}
 is a closed minimal submanifold in the sphere $\S^{p+1}(\sqrt{2})$. Denote its second fundamental form by $\td{A}$ , then
\begin{equation*}
    |A|^2=|\td{A}|^2+1,
\end{equation*}
where $A$ is the second fundamental form of $x:M^2\ra\R^{2+p}$. Then the condition \eqref{dim2-1} is equivalent to
\begin{equation}
    \f 23\leq |\td{A}|^2\leq \f 56.
\end{equation}
A well-known theorem (cf. Theorem C in \cite{KS}) implies either

\vskip 1mm

(i) $|\td{A}|^2\equiv \f 23$ and $\td{M}^2=\S^2(\sqrt{6})\hookrightarrow\S^{4}(\sqrt{2})$ is a Veronese surface.

(ii) $|\td{A}|^2\equiv \f 56$ and $\td{M}^2=\S^2(\sqrt{12})\hookrightarrow\S^6(\sqrt{2})$ is a canonical immersion.

\noindent Therefore in terms of $|A|^2$, we have two possibilities:

\vskip 1mm

(i) $|A|^2\equiv \f 53$ and the self-shrinker $x:M^2=\S^2(\sqrt{6})\ra\S^{4}(\sqrt{2})\hookrightarrow\R^{5}$ is a Veronese surface.

(ii) $|A|^2\equiv \f {11}6$ and the self-shrinker  $x:M^2=\S^2(\sqrt{12})\ra\S^6(\sqrt{2})\hookrightarrow\R^{5}$ is a canonical immersion.
\end{proof}

\subsection{Self-shrinkers with flat normal bundle}

\begin{proof}[Proof of Theorem \ref{thm1-6}]
From the equation \eqref{a-1} in the Appendix, the condition ``flat normal bundle", i.e., $R^{\perp}=0$ and $\sigma_{\alpha\beta}=\f 1p|A|^2\delta_{\alpha\beta}$ imply
\begin{equation}\label{7-1}
    \f 12\mathcal{L}|A|^2=|\nabla A|^2+\f 1p|A|^2(p-|A|^2).
\end{equation}

Since $M$ has bounded $|A|^2$ and polynomial volume growth, from Proposition \ref{app-prop2}  in the Appendix, we have
\begin{equation}\label{7-2}
    \int_M|\nabla A|^2e^{-\f 12|x|^2}<+\infty.
\end{equation}

By \eqref{7-1} and \eqref{7-2}, $\mathcal{L}|A|^2$ has finite weighted integral
\begin{equation}
    \int_M(\mathcal{L}|A|^2)e^{-\f 12|x|^2}~<~+\infty,
\end{equation}
then the self-adjointness of the operator $\mathcal{L}$ (Corollary \ref{cor2-1}) implies
\begin{align*}
    0=&\f 12\int_M(\mathcal{L}|A|^2)e^{-\f 12|x|^2}\\
    =&\int_M\left(|\nabla A|^2+\f 1p|A|^2(p-|A|^2)\right)e^{-\f 12|x|^2}.
\end{align*}
That is
\begin{equation}
    \int_M|\nabla A|^2e^{-\f 12|x|^2}=\f 1p\int_M|A|^2(|A|^2-p)e^{-\f 12|x|^2}.
\end{equation}
Therefore our assumption \eqref{in1-8} implies either  $|A|^2\equiv 0$ and $M^n$ is a plane; or $|A|^2\equiv p$ and $|\nabla A|\equiv 0$, noting our assumption $\sigma_{\alpha\beta}=\f 1p|A|^2\delta_{\alpha\beta}$, we can conclude that  $M$ is
\begin{equation*}
    x:M^n=N^{mp}\times\R^{n-mp}\hookrightarrow\R^{n+p},
\end{equation*}
where
\begin{equation*}
    N^{mp}=\S^m(\sqrt{m})\times\cdots\times\S^m(\sqrt{m})\hookrightarrow\S^{(m+1)p-1}(\sqrt{mp}).
\end{equation*}
This completes the proof of Theorem \ref{thm1-6}.
\end{proof}

\section{Further remarks on closed self-shrinkers}

\begin{proof}[Proof of Proposition \ref{prop1-8}]
The conditions (1) and (2) imply the closed self-shrinker is a  minimal submanifold in the sphere $\S^{n+p-1}(\sqrt{n})$ have been proved by Smoczyk \cite{Smo} and Cao-Li \cite{CL2} respectively. Now we only need to prove that the condition (3) can also imply the closed self-shrinker must be the minimal submanifold in the sphere.

Recall that for self-shrinker, we have the following equations.
\begin{align}
    &\f 12\Delta |x|^2=n-|x^{\perp}|^2=n-|H|^2,\label{eq8-1}\\
    &\f 12\mathcal{L}|x|^2=n-|x|^2.\label{eq8-2}
\end{align}
Since the self-shrinker is closed, we integrate \eqref{eq8-2} with weighted $e^{-\f 12|x|^2}$. By the self-adjointness the operator $\mathcal{L}$ (Lemma \ref{lem2-1}), we have
\begin{equation}\label{eq8-3}
    0=\f 12\int_M\mathcal{L}|x|^2e^{-\f 12|x|^2}=\int_M(n-|x|^2)e^{-\f 12|x|^2}.
\end{equation}
Then the condition (3) and \eqref{eq8-3} imply $|x|^2=n$. By using equation \eqref{eq8-1}, we obtain
\begin{equation*}
    |H|^2=|x^{\perp}|^2=n=|x|^2.
\end{equation*}
Therefore $H=-x$ and $M^n$ is a minimal submanifold in the sphere $\S^{n+p-1}(\sqrt{n})$.
\end{proof}

By applying the well-known theorems on the minimal submanifold in sphere by Ejiri \cite{Ejiri}, H. Li\cite{Li93}, Itoh \cite{Itoh,Itoh2} and Yau \cite{Yau}, Proposition \ref{prop1-8} implies the following three results.

\begin{thm}\label{thm8-1}
Let $M^n$ $(n\geq 3)$ be a closed self-shrinker in $\R^{n+p}$ which satisfies one of the three conditions in Proposition \ref{prop1-8}. If the Ricci curvature of $M$ satisfies $Ric\geq \f {n-2}n$, then $M^n$ must be one of the followings:
\begin{itemize}
  \item[(1)] $M^n=\S^n(\sqrt{n})$;
  \item[(2)] $M^n=\S^k(\sqrt{\f n2})\times \S^k(\sqrt{\f n2})\hookrightarrow\S^{n+1}(\sqrt{n})$, $n=2k,p=2$;
  \item[(3)] $M^4=P_c^2(\f 13)\hookrightarrow\S^7(2)$, $n=4,p=4$,
\end{itemize}
where $P_c^2(\f 13)$ denotes the complex projective space with sectional curvature smaller than $\f 13$.
\end{thm}

\begin{rem}
We remark that Ejiri's result \cite{Ejiri} holds for $n\geq 4$, which was extended by H. Li \cite{Li93} to 3-dimensional case.
\end{rem}

\begin{thm}
Let $M^n$ be a closed self-shrinker in $\R^{n+p}$ which satisfies one of the three conditions in Proposition \ref{prop1-8}. If the sectional curvature of $M$ satisfies $K\geq \f 1{2(n+1)}$, then $M^n$ must be one of the followings:
\begin{itemize}
  \item[(1)] $M^n=\S^n(\sqrt{n})$;
  \item[(2)] $M^n=\S^n(\sqrt{2(n+1)})\hookrightarrow\S^{n+p-1}(\sqrt{n})$.
\end{itemize}
\end{thm}

\begin{thm}\label{thm8-3}
Let $M^n$ be a closed self-shrinker in $\R^{n+p}$ which satisfies one of the three conditions in Proposition \ref{prop1-8}. If the sectional curvature of $M$ satisfies $K\geq \f {p-2}{2p-3}$, then $M^n$ must be one of the followings:
\begin{itemize}
  \item[(1)] $M^n=\S^k(\sqrt{k})\times\S^{n-k}(\sqrt{n-k})\hookrightarrow\S^{n+p-1}(\sqrt{n})$, $0\leq k\leq n$;
  \item[(2)] $M^2=\S^2(\sqrt{6})\hookrightarrow\S^4(\sqrt{2})$, $n=2,p=3$.
\end{itemize}
\end{thm}
\appendix

\section{Two formulas of Simons' type and the weighted integral estimates}

In this appendix, we give two formulas of Simons' type for self-shrinkers and the weighted integral estimates of the first and the second covariant derivatives of the second fundamental form of self-shrinkers,  which we used in subsection 5.3.

\begin{prop}
Let $x:M^n\ra\R^{n+p}$ be an immersed self-shrinker, then we have the following Simons' type formula,
\begin{align}
    \f 12\mathcal{L}|A|^2=&|\nabla A|^2+|A|^2-\sum_{\alpha,\beta}\sigma_{\alpha\beta}^2-|R^{\perp}|^2\label{a-1}\\
    \f 12\mathcal{L}|\nabla A|^2=&|\nabla^2A|^2+2|\nabla A|^2+6R_{\beta\alpha kl}h^{\alpha}_{ijk}h^{\beta}_{ijl}-\sigma_{\alpha\beta}h^{\alpha}_{ijk}h^{\beta}_{ijk}\label{a-2}\\
    &+6h^{\alpha}_{ijk}h^{\alpha}_{rjl}\left(h^{\beta}_{rk}h^{\beta}_{il}-h^{\beta}_{rl}h^{\beta}_{ik}\right)-3h^{\alpha}_{ijk}h^{\alpha}_{rij}h^{\beta}_{rl}h^{\beta}_{kl},\nonumber
\end{align}
where $\sigma_{\alpha\beta}=\sum\limits_{ij}h^{\alpha}_{ij}h^{\beta}_{ij}$, $R_{\beta\alpha kl}=\sum\limits_i(h^{\beta}_{ki}h^{\alpha}_{il}-h^{\beta}_{li}h^{\alpha}_{ik})$ is the curvature of the normal bundle and $|R^{\perp}|^2=\sum R_{\beta\alpha kl}^2$. In particular, for hypersurface self-shrinkers, that is $p=1$, we have
\begin{align}
    \f 12\mathcal{L}|A|^2=&|\nabla A|^2+|A|^2(1-|A|^2)\label{ap-1}\\
    \f 12\mathcal{L}|\nabla A|^2=&|\nabla^2A|^2-|\nabla A|^2(|A|^2-2)-\f{3}2|\nabla|A|^2|^2-3\Xi,\label{ap-2}
\end{align}
where $\Xi=h_{ijk}h_{ijr}h_{kl}h_{lr}-2h_{ikj}h_{jrl}h_{il}h_{kr}$.
\end{prop}
\begin{rem}
In the hypersurface case, the formulas \eqref{ap-1} and \eqref{ap-2} were derived by Colding-Minicozzi II \cite{CM2} and Ding-Xin \cite{DX2} respectively. In arbitrary codimension, the formula \eqref{a-1} was also proved by Ding-Xin \cite{DX2}.
\end{rem}
\begin{proof}
For an immersion $x:M^n\ra\R^{n+p}$, we have the following formulas of Simons' type (for example, see \cite{DX2,Li})
\begin{align*}
    \f 12\Delta|A|^2=&|\nabla A|^2+\sum_{i,k,\alpha}h^{\alpha}_{ik}H^{\alpha}_{,ik}+\sum_{i,j,k,\alpha,\beta}H^{\beta}h^{\beta}_{jk}h^{\alpha}_{ki}h^{\alpha}_{ij}\\
    &\qquad-\sum_{\alpha,\beta}\sigma_{\alpha\beta}^2-|R^{\perp}|^2.
\end{align*}
So for the self-shrinker, substituting \eqref{3-3} into the above equation, we obtain
\begin{align*}
    \f 12\mathcal{L}|A|^2=&\f 12\left(\Delta |A|^2-<x,\nabla |A|^2>\right)\\
    =&|\nabla A|^2+|A|^2-\sum_{\alpha,\beta}\sigma_{\alpha\beta}^2-|R^{\perp}|^2.
\end{align*}
This is \eqref{a-1}. Now we prove \eqref{a-2}, recall that for the covariant derivatives of the second fundamental form, we have the following  Ricci identities (see \cite {CL1,Li})
\begin{align}
    h^{\alpha}_{ijkl}-h^{\alpha}_{ijlk}=&h^{\alpha}_{rj}R_{rikl}+h^{\alpha}_{ir}R_{rjkl}+h^{\beta}_{ij}R_{\beta\alpha kl},\label{ricci-1}\\
    h^{\alpha}_{ijkls}-h^{\alpha}_{ijksl}=&h^{\alpha}_{rjk}R_{rils}+h^{\alpha}_{irk}R_{rjls}+h^{\alpha}_{ijr}R_{rkls}+h^{\beta}_{ijk}R_{\beta\alpha ls}.
\end{align}
Then
\begin{align*}
    \Delta h^{\alpha}_{ijk}=&h^{\alpha}_{ijkll}\\
    =&(h^{\alpha}_{ijlk}+h^{\alpha}_{rj}R_{rikl}+h^{\alpha}_{ir}R_{rjkl}+h^{\beta}_{ij}R_{\beta\alpha kl})_{l}\\
    =&h^{\alpha}_{lijlk}+h^{\alpha}_{rij}R_{rlkl}+h^{\alpha}_{lrj}R_{rikl}+h^{\alpha}_{lir}R_{rjkl}+h^{\beta}_{lij}R_{\beta\alpha kl}\\
    &\quad +h^{\alpha}_{rjl}R_{rikl}+h^{\alpha}_{irl}R_{rjkl}+h^{\beta}_{ijl}R_{\beta\alpha kl}\\
    &\quad+h^{\alpha}_{rj}(R_{rikl})_l+h^{\alpha}_{ir}(R_{rjkl})_l+h^{\beta}_{ij}(R_{\beta\alpha kl})_l\\
    =&H^{\alpha}_{,ijk}+h^{\alpha}_{rik}R_{rj}+h^{\alpha}_{rij}R_{rk}+h^{\alpha}_{lrk}R_{rijl}+h^{\beta}_{lik}R_{\beta\alpha jl}\\
    &\quad +2h^{\alpha}_{rjl}R_{rikl}+2h^{\alpha}_{ril}R_{rjkl}+2h^{\beta}_{lij}R_{\beta\alpha kl}\\
    &\quad+h^{\alpha}_{rj}(R_{rikl})_l+h^{\alpha}_{ir}(R_{rjkl})_l+h^{\beta}_{ij}(R_{\beta\alpha kl})_l\\
    &\quad+h^{\alpha}_{ri}(R_{rljl})_k+h^{\alpha}_{rl}(R_{rijl})_k+h^{\beta}_{li}(R_{\beta\alpha jl})_k,
\end{align*}
where we used in last equality
\begin{align*}
 h^\alpha_{lijlk}=&[h^\alpha_{lilj}+h^\alpha_{ri}R_{rljl}+h^\alpha_{lr}R_{rijl}+h^\beta_{li}R_{\beta\alpha jl}]_k\\
 =&H^\alpha_{,ijk}+[h^\alpha_{rik}R_{rj}+h^\alpha_{lrk}R_{rijl}+h^\alpha_{lik}R_{\beta\alpha jl}]\\
 &\qquad h^\alpha_{ri}(R_{rj})_k+h^\alpha_{lr}(R_{rijl})_k+h^\beta_{li}(R_{\beta\alpha jl})_k.
\end{align*}
For self-shrinkers, we have by use of \eqref{3-1}, \eqref{3-2} and \eqref{3-3}
\begin{align}
    H^{\alpha}_{,ijk}=&h^{\alpha}_{ijlk}<x,e_l>+2h^{\alpha}_{ijk}-H^{\beta}_{,k}h^{\alpha}_{il}h^{\beta}_{lj}\nonumber\\
    &\quad -H^{\beta}\left(h^{\alpha}_{ijl}h^{\beta}_{kl}+h^{\alpha}_{ikl}h^{\beta}_{jl}+h^{\alpha}_{il}h^{\beta}_{ljk}\right)
\end{align}
Then we obtain
\begin{align}
    \f 12\mathcal{L}|\nabla A|^2=&\f 12(\Delta |\nabla A|^2-<x,\nabla |\nabla A|^2>)\nonumber\\
    =&|\nabla^2A|^2+h^{\alpha}_{ijk}\Delta h^{\alpha}_{ijk}-<x,e_l>h^{\alpha}_{ijk}h^{\alpha}_{ijkl}\nonumber\\
    =&|\nabla^2A|^2+2|\nabla A|^2+<x,e_l>h^{\alpha}_{ijk}(h^{\alpha}_{ijlk}-h^{\alpha}_{ijkl})\nonumber\\
    &\quad -h^{\alpha}_{ijk}H^{\beta}_{,k}h^{\alpha}_{il}h^{\beta}_{lj}-h^{\alpha}_{ijk}H^{\beta}\left(h^{\alpha}_{ijl}h^{\beta}_{kl}+h^{\alpha}_{ikl}h^{\beta}_{jl}+h^{\alpha}_{il}h^{\beta}_{ljk}\right)\nonumber\\
    &\quad+h^{\alpha}_{ijk}\left(h^{\alpha}_{rik}R_{rj}+h^{\alpha}_{rij}R_{rk}+h^{\alpha}_{lrk}R_{rijl}+h^{\beta}_{lik}R_{\beta\alpha jl}\right.\label{a-5}\\
    &\quad +2h^{\alpha}_{rjl}R_{rikl}+2h^{\alpha}_{ril}R_{rjkl}+2h^{\beta}_{lij}R_{\beta\alpha kl}\nonumber\\
    &\quad+h^{\alpha}_{rj}(R_{rikl})_l+h^{\alpha}_{ir}(R_{rjkl})_l+h^{\beta}_{ij}(R_{\beta\alpha kl})_l\nonumber\\
    &\quad\left.+h^{\alpha}_{ri}(R_{rljl})_k+h^{\alpha}_{rl}(R_{rijl})_k+h^{\beta}_{li}(R_{\beta\alpha jl})_k\right)\nonumber
\end{align}
By the Ricci identity \eqref{ricci-1}, Gauss equation \eqref{gauss-1} and the equation \eqref{3-2}, a direct calculation to check
\begin{align}
    <x,e_l>h^{\alpha}_{ijk}(h^{\alpha}_{ijlk}-h^{\alpha}_{ijkl})=&3h^{\alpha}_{ijk}h^{\alpha}_{rj}h^{\beta}_{ik}H^{\beta}_{,r}-2h^{\alpha}_{ijk}h^{\alpha}_{rj}h^{\beta}_{rk}H^{\beta}_{,i}\nonumber\\
    &\quad-h^{\alpha}_{ijk}h^{\beta}_{ij}h^{\beta}_{rk}H^{\alpha}_{,r}\label{a-6}
\end{align}
and the last four lines of \eqref{a-5} is equal to
\begin{align}
    &h^{\alpha}_{ijk}\left(6h^{\alpha}_{rjl}h^{\beta}_{rk}h^{\beta}_{il}-6h^{\alpha}_{rjl}h^{\beta}_{rl}h^{\beta}_{ik}+2h^{\alpha}_{rij}h^{\beta}_{rk}H^{\beta}-3h^{\alpha}_{rij}h^{\beta}_{rl}h^{\beta}_{kl}\right.\nonumber\\
    &\quad + 6h^{\beta}_{lij}h^{\beta}_{kr}h^{\alpha}_{lr}-6h^{\beta}_{lij}h^{\beta}_{lr}h^{\alpha}_{kr}+3h^{\alpha}_{ri}h^{\beta}_{rj}H^{\beta}_{,k}-3h^{\alpha}_{rj}h^{\beta}_{ki}H^{\beta}_{,r}\label{a-7}\\
    &\quad \left.+h^{\alpha}_{ri}h^{\beta}_{rjk}H^{\beta}-h^{\alpha}_{rl}h^{\beta}_{rl}h^{\beta}_{ijk}+h^{\beta}_{ij}h^{\beta}_{kr}H^{\alpha}_{,r}\right)\nonumber
\end{align}
Put \eqref{a-6} and \eqref{a-7} into \eqref{a-5}, we get the formula \eqref{a-2}.
\end{proof}

\begin{prop}\label{app-prop2}
Let $x:M^n\ra\R^{n+p}$ be a complete immersed self-shrinker and with polynomial volume growth, if $|A|^2$ is bounded on $M$, then
\begin{align}
    &\int_M|\nabla A|^2e^{-\f 12|x|^2}<+\infty\label{a-3},\\
    &\int_M|\nabla^2 A|^2e^{-\f 12|x|^2}<+\infty\label{a-4}.
\end{align}
\end{prop}
\begin{rem}
For the hypersurface case, this weighted integral estimate was proved by Ding-Xin (\cite{DX2}).
\end{rem}
\begin{proof}
Let $\eta$ be a cut-off function with compact support on $M$, by the self-adjointness of $\mathcal{L}$ and \eqref{a-1}, we have
\begin{align*}
    \int_M|\nabla A|^2\eta^2e^{-\f 12|x|^2}=&\int_M\left(\sum_{\alpha,\beta}\sigma_{\alpha\beta}^2+|R^{\perp}|^2-|A|^2+\f 12\mathcal{L}|A|^2\right)\eta^2e^{-\f 12|x|^2}\\
    \leq&C\int_M|A|^4\eta^2e^{-\f 12|x|^2}-\f 12\int_M\nabla\eta^2\nabla|A|^2e^{-\f 12|x|^2}\\
    \leq&C\int_M|A|^4\eta^2e^{-\f 12|x|^2}+\f 12 \int_M|\nabla A|^2\eta^2e^{-\f 12|x|^2}\\
    &\qquad +2\int_M|A|^2|\nabla\eta|^2e^{-\f 12|x|^2}.
\end{align*}
Then
\begin{equation*}
    \int_M|\nabla A|^2\eta^2e^{-\f 12|x|^2}\leq 2C\int_M|A|^4\eta^2e^{-\f 12|x|^2}+4\int_M|A|^2|\nabla\eta|^2e^{-\f 12|x|^2}.
\end{equation*}
Since $M$ has polynomial volume growth and bounded $|A|^2$, by the dominated convergence theorem, the above inequality implies
\begin{equation*}
    \int_M|\nabla A|^2e^{-\f 12|x|^2}<+\infty.
\end{equation*}

Using the formula \eqref{a-2} and a similar argument as above, we get
\begin{align*}
   \int_M|\nabla^2 A|^2\eta^2e^{-\f 12|x|^2}\leq &C\int_M|A|^2|\nabla A|^2e^{-\f 12|x|^2}+4\int_M|\nabla A|^2|\nabla \eta|^2e^{-\f 12|x|^2}.
\end{align*}
Then by the boundness of $|A|^2$, \eqref{a-3}, the dominated convergence theorem and the above equality imply
\begin{equation*}
    \int_M|\nabla^2 A|^2e^{-\f 12|x|^2}<+\infty.
\end{equation*}
\end{proof}

\bibliographystyle{Plain}

\begin{thebibliography}{10}

\bibitem{AbL} U. Abresch and J. Langer, {\it The normalized curve shortening flow and homothetic solutions}, J. Diff. Geom., {\bf 23} (1986),no. 2, 175-196.

\bibitem{ALW} B. Andrews, H. Li and Y. Wei, {\it $\mathcal{F}$-stability for self-shrinking solutions to mean curvature flow}, arXiv:1204.5010v1.

\bibitem{CL1} L. F. Cao and H. Li, {\it $r$-minimal submanifolds in space forms}, Ann. Global Anal. Geom.,
{\bf 32} (2007), 311-341.

\bibitem{CL2} H. -D. Cao and H. Li, {\it A gap theorem for self-shrinkers of the mean curvature flow in arbitrary
codimension}, to appear in Calc. Var., (arXiv: 1101.0516v1).


\bibitem{Cheng-P} Q.-M. Cheng and Y. Peng, {\it Complete self-shrinkers of the mean curvature flow}, arXiv: 1202.1053v1.

\bibitem{ChZ} X. Cheng and D. Zhou, {\it Volume estimate about shrinkers}, arXiv: 1106.4950v1

\bibitem{CM2} T. H. Colding and W.P. Minicozzi II, {\it Generic mean curvature flow
I: generic singularities}, Ann. of Math., {\bf 175}(2) (2012), 755-833.


\bibitem{DX1} Q. Ding and Y. L. Xin, {\it On Chern's problem for rigidity of minimal hypersurfaces in the spheres}, Adv. Math., \textbf{227}(2011),131-145.

\bibitem{DX2} Q. Ding and Y. L. Xin, {\it The rigidity theorems of self-shrinkers}, arXiv:1105:4962v1.

\bibitem{Ejiri} N. Ejiri, {\it Compact minimal submanifolds of a sphere
with positive Ricci curvature}, J. Math. Soc. Japan, \textbf{31}(1979)(251-256).

\bibitem{H2} G. Huisken, {\it Asymptotic behavior for singularities of the mean
curvature flow,} J. Differential Geom., {\bf 31} (1990), no. 1, 285-299.

\bibitem{H3} G. Huisken, {\it Local and global behaviour of hypersurfaces moving by
mean curvature,} Proc. Sympos. Pure Math.,  {\bf 54}
(1993), Amer. Math. Soc.

\bibitem{I} T. Ilmanen, {\it Singularities of mean curvature flow of surfaces}, preprint, 1995, \url{http://www.math.ethz.ch/~ilmanen/papers/pub.html}.

\bibitem{Itoh} T. Itoh, {\it On Veronese manifolds}, J. Math. Soc. Japan, \textbf{27}(1975),497-506.

\bibitem{Itoh2} T. Itoh, {\it Addendum to my paper ``On Veronese manifolds''}, J. Math. Soc. Japan, \textbf{30}(1978),73-74.

\bibitem{KM} S. Kleene and N.M. M{\o}ller, {\it Self-shrinkers with a rotaion symmetry,} arXiv:1008.1609v1.

\bibitem{KS} M. Kozlowski and U. Simon, {\it Minimal immersion of 2-manifolds into spheres}, Math. Z., \textbf{186}(1984),377-382.

\bibitem{LeS} Nam Q. Le and N. Sesum, {\it Blow-up rate of the mean curvature during the mean curvature flow and a gap theorem for self-shrinkers}, Comm. Anal. Geom., \textbf{19}(2011), no.4, 633-660.

\bibitem{Li93} H. Li, {\it Curvature pinching for odd-dimensional minimal submanifolds in a sphere}, Publ. Inst. Math. (Beograd), \textbf{53}(1993), 122-132.

\bibitem{Li} H. Li, {\it Willmore submanifolds in a sphere}, Math. Res. Lett.,  {\bf 9} (2002), 771-790.

\bibitem{LS} H. Li and U. Simon, {\it Quantization of curvature for compact surfaces in $S^n$}, Math. Z., \textbf{245}(2003), 201-216.

\bibitem{LW} H. Li and Y. Wei, {\it Lower volume growth estimates for self-shrinkers of mean curvature flow,} arXiv: 1112.0828, to appear in Proceedings of 
 Amer. Math. Soc.

\bibitem{SSY} R. Schoen, L. Simon and S.-T. Yau, {\it Curvature estimates for minimal hypersurfaces},  Acta Math., {\bf 134} (1975), 275-288.

\bibitem{Si} J. Simons, {\it Minimal varieties in Riemannian manifolds}, Ann. of Math., \textbf{88}(1968), no.1,62-105.

\bibitem{Smo} K. Smoczyk, {\it Self-shrinkers of the mean curvature flow in arbitrary codimension,} Int. Math. Res. Not., {\bf 48} (2005), 2983-3004.

\bibitem{White} B. White, {\it Evolution of curves and surfaces by mean curvature}, ICM 2002, Vol. I.  525-538.

\bibitem{Yau}S. -T. Yau, {\it Submanifolds with constant mean curvature I, II}, Amer.  J. Math.,
{\bf 96} (1974), 346-366., {\bf 97}(1975),76-100.

\bibitem{Zhu} X. P. Zhu, {\it Lectures on mean curvature flows}, Studies in Adv. Math, AMS and IP,2002.

\end{thebibliography}

\end{document}